\documentclass[12pt]{article}
\usepackage[latin1]{inputenc}
\usepackage[british]{babel}
\usepackage{cmap}
\usepackage{lmodern}

\usepackage{amssymb, amsmath, amsthm}
\usepackage[a4paper,top=25mm,bottom=25mm,left=25mm,right=25mm]{geometry}
\usepackage{ragged2e}

\usepackage{authblk} 
\usepackage{pifont}
\usepackage{graphicx}
\usepackage[usenames,dvipsnames,svgnames,table]{xcolor}
\usepackage[figuresright]{rotating}
\usepackage{xtab} 
\usepackage{longtable} 
\usepackage{multirow}
\usepackage{footnote}
\usepackage[stable]{footmisc}
\usepackage{chngpage} 
\usepackage{pdflscape} 
\usepackage[nottoc,notlot,notlof]{tocbibind} 

\usepackage{pgfplots}
\pgfplotsset{every tick label/.append style={font=\footnotesize}}
\pgfplotsset{compat=1.18}
\usepackage{setspace}

\usepackage{array}
\newcolumntype{K}[1]{>{\centering\arraybackslash$}p{#1}<{$}}

\makesavenoteenv{tabular}
\usepackage{tabularx}
\usepackage{booktabs}
\usepackage{threeparttable}
\usepackage[referable]{threeparttablex} 
\newcolumntype{R}{>{\raggedleft\arraybackslash}X}
\newcolumntype{L}{>{\raggedright\arraybackslash}X}
\newcolumntype{C}{>{\centering\arraybackslash}X}
\newcolumntype{A}{>{\columncolor{gray!25}}C}
\newcolumntype{a}{>{\columncolor{gray!25}}c}

\newlength{\tablen}

\usepackage{dcolumn} 
\newcolumntype{.}{D{.}{.}{-1}}

\usepackage{tikz}
\usetikzlibrary{arrows, arrows.meta, calc, matrix, patterns, positioning, trees}
\tikzset{vertex/.style={circle,draw}, edge/.style={->,> = latex'}}
\usepackage[semicolon]{natbib}
\usepackage[hyphens]{url}
\usepackage{hyperref} 
\hypersetup{
  colorlinks   = true,    		
  urlcolor     = blue,    		
  linkcolor    = blue,    		
  citecolor    = ForestGreen	
}
\usepackage{microtype}
\usepackage[justification=centering]{caption} 

\usepackage[labelformat=simple]{subcaption}

\DeclareCaptionLabelFormat{parenthesis}{(#2)}
\makeatletter
\renewcommand\p@subfigure{\arabic{figure}.}
\makeatother

\DeclareCaptionLabelFormat{parenthesis}{(#2)}
\makeatletter
\renewcommand\p@subtable{\arabic{table}.}
\makeatother

\usepackage{enumitem}

\setlist[itemize]{leftmargin=2.5\parindent}
\setlist[enumerate]{leftmargin=2.5\parindent}

\newenvironment{customlegend}[1][]{%
	\begingroup
	\csname pgfplots@init@cleared@structures\endcsname
	\pgfplotsset{#1}%
    }{%
	\csname pgfplots@createlegend\endcsname
	\endgroup
    }%
\def\addlegendimage{\csname pgfplots@addlegendimage\endcsname}

\theoremstyle{plain}

\theoremstyle{definition}

\newtheorem{example}{Example}

\theoremstyle{remark}

\makeatletter
\let\@fnsymbol\@alph
\makeatother

\def\keywords{\vspace{.5em} 
{\noindent \textit{Keywords}: }}

\def\AMS{\vspace{.5em} 
{\noindent \textbf{\emph{MSC} class}: }}

\def\JEL{\vspace{.5em} 
{\noindent \textbf{\emph{JEL} classification number}: }}

\title{Right-left asymmetry of the eigenvector method: \\ A simulation study}
\author{\href{https://sites.google.com/view/laszlocsato}{L\'aszl\'o Csat\'o}\thanks{~E-mail: \emph{laszlo.csato@sztaki.hu}} }
\affil{HUN-REN Institute for Computer Science and Control (HUN-REN SZTAKI) \\
Laboratory on Engineering and Management Intelligence \\
Research Group of Operations Research and Decision Systems}
\affil{Corvinus University of Budapest (BCE) \\
Institute of Operations and Decision Sciences \\
Department of Operations Research and Actuarial Sciences}
\affil{Budapest, Hungary}

\def\Dedication{
\begin{small}
{\noindent
``\emph{A suggestion for the prioritization of alternatives using the Perron-Frobenius right eigenvector of a pairwise comparison matrix has recently been made by T.~Saaty. We note that use of the left eigenvectors is equally justified (as long as order is reversed).}''\footnote{~Source: \citet[p.~61]{JohnsonBeineWang1979}.}
}
\end{small}

\vspace{0.5cm} 
\justify }

\begin{document}

\maketitle
\thispagestyle{empty}
\Dedication

\begin{abstract}
\noindent
The eigenvalue method, suggested by the developer of the extensively used Analytic Hierarchy Process methodology, exhibits right-left asymmetry: the priorities derived from the right eigenvector do not necessarily coincide with the priorities derived from the reciprocal left eigenvector. This paper offers a comprehensive numerical experiment to compare the two eigenvector-based weighting procedures and their reasonable alternative of the row geometric mean with respect to four measures. The underlying pairwise comparison matrices are constructed randomly with different dimensions and levels of inconsistency. The disagreement between the two eigenvectors turns out to be not always a monotonic function of these important characteristics of the matrix. The ranking contradictions can affect alternatives with relatively distant priorities. The row geometric mean is found to be almost at the midpoint between the right and inverse left eigenvectors, making it a straightforward compromise between them.


\AMS{90-10, 90B50, 91B08}

\JEL{C44, D71}
\end{abstract}

\clearpage

\section{Introduction} \label{Sec1}

The Analytic Hierarchy Process (AHP) is one of the most popular decision-making techniques since it has been introduced by \citet{Saaty1977, Saaty1980}. It has a number of successful applications \citep{BhushanRai2007, FormanGass2001, VaidyaKumar2006, Vargas1990} and, simultaneously, several flaws identified in the literature \citep{Csato2017c, CsatoPetroczy2021, GenestLapointeDrury1993, MunierHontoria2021, PetroczyCsato2021}.

The current paper deals with an issue in priority derivation from a given pairwise comparison matrix. Saaty has suggested using the right eigenvector for this purpose but there are many other methods \citep{ChooWedley2004}. In particular, \citet{JohnsonBeineWang1979} argue for the componentwise reciprocal of the left eigenvector as written in the motto above.
Another strong competitor is the logarithmic least squares or row geometric mean \citep{CrawfordWilliams1985}, mainly due to its strong axiomatic foundations \citep{Fichtner1984, Fichtner1986, BarzilaiCookGolany1987, Barzilai1997, LundySirajGreco2017, Csato2018b, BozokiTsyganok2019, Csato2019a}.

In the following, these solutions will be compared using a Monte Carlo simulation approach, that is, the priority vectors of the three weighting procedures are evaluated on the basis of a large set of random pairwise comparison matrices. Even though a similar exercise has been attempted at least twice in the previous literature \citep{BozokiRapcsak2008, IshizakaLusti2006}, independently of each other, we make important contributions compared to both preliminary studies:
\begin{itemize}
\item
\citet{IshizakaLusti2006} limit the investigation to at most seven alternatives and matrices with an acceptable level of inconsistency. The simulation has some drawbacks because
(1) the number of matrices for a given order is 500, which might be insufficient to derive robust results; and
(2) the matrix entries can only be integers.
Last but not least, the authors focus exclusively on the number of ranking contradictions and do not consider the differences in the priorities derived from the three methods.
\item
\citet{BozokiRapcsak2008} examine the case of five alternatives. All matrix entries are chosen randomly, hence, the number of matrices with an acceptable level of inconsistency remains rather low. The authors only consider the frequency of rank reversals.
\end{itemize}
Here, 3 million pairwise comparison matrices are generated for any given number of alternatives between four and nine such that each interval of Saaty's inconsistency ratio \citep{Saaty1977} with a length of $0.005$ contains at least one thousand instances. The right and inverse left eigenvectors, as well as the row geometric mean, are evaluated according to four measures: the Euclidean and Chebyshev distances, as well as the Kendall rank correlation coefficients of the normalised weight vectors, and the maximal ratios of the priorities corresponding to one alternative. Among them, only Kendall rank correlation depends on the ranking of the alternatives.

Our comprehensive analysis yields several interesting results, some of them at least partially contradicting previous findings:
\begin{itemize}
\item
The row geometric mean is revealed to be an excellent compromise between the right and inverse left eigenvectors as it is almost at the midpoint between them, especially at a low level of inconsistency;
\item
The differences between the three priority deriving methods do not always increase with the level of inconsistency if the latter is relatively high (we refine the conclusion of \citet[p.~398]{IshizakaLusti2006});
\item
The differences between the three priority deriving methods do not always increase with the number of alternatives (we refine the conclusion of \citet[p.~398]{IshizakaLusti2006});
\item
Three examples illustrate that
(1) rank reversal between the right and inverse left eigenvectors may emerge for a slightly perturbed consistent matrix;
(2) the right and inverse left eigenvectors can lead to a fully reversed order of the alternatives;
(3) rank reversal between the right and inverse left eigenvectors might occur even if the priorities of two alternatives are distant (this denies the conclusion of \citet[p.~398]{IshizakaLusti2006}).
\end{itemize}

The main reason for the different conclusions compared to \citet{IshizakaLusti2006} resides in the extension of our analysis to
(a) a wider interval of inconsistency (inconsistency ratio $\mathit{CR} < 0.5$ rather than $\mathit{CR} < 0.1$);
(b) a higher number of alternatives (up to nine instead of seven);
(c) a broader set of comparison metrics.
The expansion of the range of inconsistency can be justified since $\mathit{CR} < 0.1$ is an inflexible criterion and is too restrictive when the size of the matrix increases \citep{AlonsoLamata2006}.
Furthermore, the results are based on a much higher number of random matrices.

The remainder of the study is organised as follows.
The mathematical background is presented in Section~\ref{Sec2}. Section~\ref{Sec3} reviews the related literature, and Section~\ref{Sec4} outlines the simulation experiment. Section~\ref{Sec5} contains the main results. The paper is finished with a concise discussion in Section~\ref{Sec6}.

\section{Theoretical background} \label{Sec2}

An $n \times n$ matrix $\mathbf{A} = \left[ a_{ij} \right]$ is called a pairwise comparison matrix if it is positive ($a_{ij} > 0$ for all $i,j$) and reciprocal ($a_{ji} = 1 / a_{ij}$ for all $i,j$). Its entry $a_{ij}$ quantifies how many times alternative $i$ is better/more important compared to alternative $j$.
An important property of a pairwise comparison matrix is \emph{consistency}: it is called consistent if $a_{ik} = a_{ij} a_{jk}$ for all $i,j,k$; otherwise, it is called inconsistent.

Pairwise comparison matrices are mainly used to derive priorities for the alternatives. In the case of a consistent matrix, this is almost trivial since the matrix is generated by an appropriate weight vector, that is, there exists a vector $\mathbf{w} = \left[ w_i \right]$ such that $a_{ij} = w_i / w_j$ for all $i,j$. For inconsistent matrices, several weighting techniques have been proposed in the literature \citep{ChooWedley2004}.
Probably the two most popular procedures are the logarithmic least squares (LLSM) or row geometric mean \citep{CrawfordWilliams1985, DeGraan1980, deJong1984, Rabinowitz1976, WilliamsCrawford1980} and the eigenvector \citep{Saaty1977} methods.

The logarithmic least squares method minimises the aggregated distances of the approximations in a logarithmic sense:
\[
\sum_{i=1}^n \sum_{j=1}^n \left[ \log \left( a_{ij} \right) - \log \left( \frac{w_i}{w_j} \right) \right]^2 \to \min.
\]
The solution is provided by the geometric means of row elements, namely:
\[
\frac{w_i}{w_j} = \frac{\sqrt[n]{\prod_{k=1}^n a_{ik}}}{\sqrt[n]{\prod_{k=1}^n a_{jk}}}.
\]
The corresponding weight vector is denoted by $\mathbf{w}^{\mathit{RGM}}$.

The eigenvector method is based on the right eigenvector associated with the dominant eigenvalue $\lambda_{\max}$ of the pairwise comparison matrix:
\[
\mathbf{A} \mathbf{w}^R = \lambda_{\max} \mathbf{w}^R.
\]

However, the matrix has also a left eigenvector associated with the same dominant eigenvalue $\lambda_{\max}$:
\[
\mathbf{w}^L \mathbf{A} = \lambda_{\max} \mathbf{w}^L.
\]

If the pairwise comparison matrix is consistent, the componentwise inverse of the left eigenvector is the right eigenvector: $w_i^L = 1 / w_i^R$. Therefore, it is reasonable to use the componentwise inverse of the left eigenvector, denoted by $\mathbf{w}^{-L}$ in the following, to derive the priorities \citep{JohnsonBeineWang1979}.

Finally, the (geometric) mean of the right and inverse left eigenvectors $\mathbf{w}^{\mathit{RL}} = \left[ w_i^{\mathit{RL}} \right]$ can be defined as
\[
w_i^{\mathit{RL}} = w_i^R w_i^{-L}
\]
for all $1 \leq i \leq n$.

Several inconsistency measures have been suggested in the literature \citep{Brunelli2018}. In this paper, we use the first index proposed by Saaty \citep{Saaty1977}:
\[
\mathit{CI} = \frac{\lambda_{\max}-n}{n-1}.
\]
The value of $\mathit{CI}$ is compared to the average $\mathit{CI}$ of a high number of randomly generated pairwise comparison matrices, which is denoted by $\mathit{RI}$, in order to get the inconsistency ratio $\mathit{CR} = \mathit{CI} / \mathit{RI}$.
According to Saaty, a pairwise comparison matrix can be accepted if $\mathit{CR}$ does not exceed $0.1$. A statistical interpretation of the 10\% rule is provided by \citet{Vargas1982}.

\section{The significance of the problem} \label{Sec3}

A conceptual weakness of the eigenvector method is the issue of asymmetry \citep{BozokiRapcsak2008}. The entry $a_{ij}$ of a pairwise comparison matrix $\mathbf{A}$ gives the numerical answer to the question ``How much does alternative $i$ dominate alternative $j$?''. However, one can equivalently ask the reciprocal question of ``How much does alternative $j$ dominate alternative $i$?'' to arrive at matrix $\mathbf{A}^\top$. The latter approach produces the right eigenvector of $\mathbf{A}^\top$, which is the elementwise reciprocal of the right eigenvector of $\mathbf{A}$.

A crucial motivation behind the eigenvector method comes from the property of consistent pairwise comparison matrices that $a_{ij} = w_{i}/w_{j}$ and $w_{i} = a_{ij} w_{j}$ for all $i,j$. This system of linear equations leads almost directly to the matrix equation $\lambda \mathbf{w} = \mathbf{A} \mathbf{w}$. However, $a_{ij} = w_{i}/w_{j}$ can be equivalently written as $w_{j} = a_{ji} w_{i}$ for all $i,j$, which implies the matrix equation $\mathbf{A}^\top \mathbf{w} = \lambda \mathbf{w}$. Obviously, the right eigenvector of $\mathbf{A}^\top$ is the left eigenvector of $\mathbf{A}$. However, the meaning of the entries in $\mathbf{A}^\top$ is the opposite compared to the entries of $\mathbf{A}$.
Consequently, ``\emph{we might just as well use a left eigenvector to prioritize in the general case as long as order reversal is allowed for}'' \citep[p.~62]{JohnsonBeineWang1979}.

It might also happen that the decision-maker misinterprets the task, and provides all pairwise comparisons in a reversed order \citep{DoddDoneganMcMaster1995}. This reversal results in the emergence of the left rather than the right eigenvector.
Hence, ``\emph{if these two vectors are not componentwise mutual inversions, it is impossible to say which is 'correct'} \citep[p.~88]{DoddDoneganMcMaster1995}.''

In order to highlight further why the right-left asymmetry can be important, a simple group decision-making problem with two sets of preferences is considered.

\begin{figure}[t!]
\centering

\begin{subfigure}{0.4\textwidth}
\centering
\caption{The preferences of DM1}
\label{Fig1a}
\begin{tikzpicture}[scale=1, auto=center, transform shape, >=triangle 45]
  \node[vertex] (n1) at (135:2) {S};
  \node[vertex] (n2) at (45:2)  {T};
  \node[vertex] (n3) at (315:2) {U};
  \node[vertex] (n4) at (225:2) {V};

  \draw[edge] [-{Latex[length=3mm]}] (n1) -- (n2) node[midway, above] {$1$};
  \draw[edge] [-{Latex[length=3mm]}] (n1) -- (n3) node[near start, above] {$1$};
  \draw[edge] [-{Latex[length=3mm]}] (n1) -- (n4) node[midway, left]  {$9$};
  \draw[edge] [-{Latex[length=3mm]}] (n2) -- (n3) node[midway, right] {$2$};
  \draw[edge] [-{Latex[length=3mm]}] (n2) -- (n4) node[near start, above] {$5$};
  \draw[edge] [-{Latex[length=3mm]}] (n3) -- (n4) node[midway, below] {$9$};
\end{tikzpicture}
\end{subfigure}
\begin{subfigure}{0.4\textwidth}
\centering
\caption{The preferences of DM2}
\label{Fig1b}
\begin{tikzpicture}[scale=1, auto=center, transform shape, >=triangle 45]
  \node[vertex] (n1) at (135:2) {S};
  \node[vertex] (n2) at (45:2)  {T};
  \node[vertex] (n3) at (315:2) {U};
  \node[vertex] (n4) at (225:2) {V};

  \draw[edge] [-{Latex[length=3mm]}] (n2) -- (n1) node[midway, above] {$1$};
  \draw[edge] [-{Latex[length=3mm]}] (n3) -- (n1) node[near start, above] {$1$};
  \draw[edge] [-{Latex[length=3mm]}] (n4) -- (n1) node[midway, left]  {$9$};
  \draw[edge] [-{Latex[length=3mm]}] (n3) -- (n2) node[midway, right] {$2$};
  \draw[edge] [-{Latex[length=3mm]}] (n4) -- (n2) node[near start, above] {$5$};
  \draw[edge] [-{Latex[length=3mm]}] (n4) -- (n3) node[midway, below] {$9$};
\end{tikzpicture}
\end{subfigure}

\caption{The preferences of two decision-makers in Example~\ref{Examp1}}
\label{Fig1}
\end{figure}
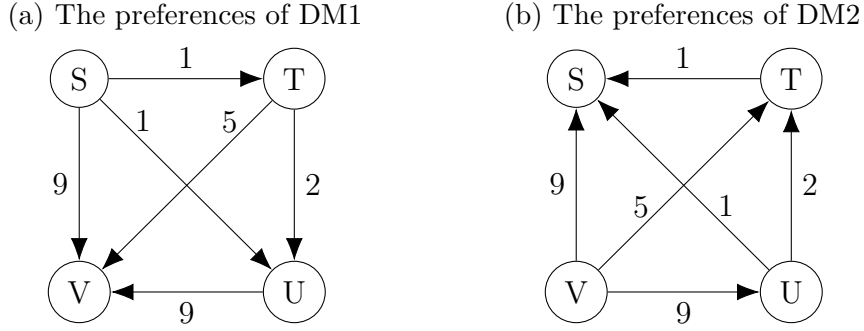

\begin{example} \label{Examp1}
Figure~\ref{Fig1} presents the numerical preferences of two decision-makers (DMs) for four alternatives. For instance, alternative S is judged to be twice more important as alternative T by DM1, but S is only half more important compared to T according to DM2. Indeed, the preferences of DM1 and DM2 are exactly the opposite. Therefore, it is natural to assign the same priorities to all alternatives S, T, U, V on the basis of the aggregated preferences of the two DMs.

Take the associated pairwise comparison matrices $\mathbf{A}$ and $\mathbf{B}$ of DM1 and DM2, respectively:
\[
\mathbf{A} = \left[
\begin{array}{K{2em} K{2em} K{2em} K{2em}}
    1     & 1     & 1     & 9     \\
    1     & 1     & 2     & 5     \\
    1     & 1/2   & 1     & 9     \\
    1/9   & 1/5   & 1/9   & 1     \\
\end{array}
\right]
\qquad \qquad \qquad
\mathbf{B} = \left[
\begin{array}{K{2em} K{2em} K{2em} K{2em}}
    1     & 1     & 1     & 1/9 \\
    1     & 1     & 1/2   & 1/5 \\
    1     & 2     & 1     & 1/9 \\
    9     & 5     & 9     & 1     \\
\end{array}
\right].
\]
The corresponding right eigenvectors are as follows (the sum of priorities are normalised to 100):
\[
\mathbf{w}^R \left( \mathbf{A} \right) = \left[
\begin{array}{K{2.5em} K{2.5em} K{2.5em} K{2.5em}}
    32.42 & 35.02 & 28.21 & 4.35 \\
\end{array}
\right];
\]
\[
\mathbf{w}^R \left( \mathbf{B} \right) = \left[
\begin{array}{K{2.5em} K{2.5em} K{2.5em} K{2.5em}}
    8.86  & 9.05  & 11.04 & 71.05 \\
\end{array}
\right].
\]

In group decision-making, there are essentially two ways to derive aggregated priorities for the alternatives:
(1) aggregating the individual pairwise comparison matrices, and using a weighting method for the aggregated matrix \citep{AczelSaaty1983}; or
(2) using a weighting method for the individual pairwise comparison matrices, and aggregating these priorities \citep{BasakSaaty1993}.
\citet{AczelSaaty1983} provide an axiomatic approach to show that individual matrices should be aggregated by the geometric mean when the aggregation of $\mathbf{A}$ and $\mathbf{B} = \mathbf{A}^\top$ results in a consistent matrix in which all entries are one. Then each alternative should have the same priority.

On the other hand, any reasonable aggregation procedure of $\mathbf{w}^R \left( \mathbf{A} \right)$ and $\mathbf{w}^R \left( \mathbf{B} \right)$ results in a higher weight for alternative T than for alternative S since the second entry of both individual weight vectors is higher than the first entry.

Note that in this case, the disturbingly different implications of techniques (1) and (2) are caused exclusively by the difference between the right and inverse left eigenvectors.
\end{example}

According to Example~\ref{Examp1}, a potential cause of the well-known difference between aggregation procedures (1) and (2) for the eigenvector method can be the right-left asymmetry.
Thus, using the right eigenvector can be strongly debated in group decision-making if the right and inverse left eigenvectors differ.
Based on these arguments, we agree with \citet[p.~62]{JohnsonBeineWang1979} that ``\emph{there is no reason to believe that utilization of a right eigenvector yields a `better' scheme than the left}.''

\section{Literature review} \label{Sec4}

While the components in the normalised left eigenvector of a pairwise comparison matrix with three alternatives are the reciprocals of the components of the normalised right eigenvector, this is not necessarily the case when the number of alternatives is at least four. This problem has been identified first by \citet{JohnsonBeineWang1979}, who presented the following example:
\[
\mathbf{A} = \left[
\begin{array}{K{2em} K{2em} K{2em} K{2em}}
    1     	& 3		  	& 1/3   	& 1/2   \\
    1/3		& 1       	& 1/6		& 2 \\
    3		& 6			& 1      	& 1 \\
    2	 	& 1/2	  	& 1			& 1 \\
\end{array}
\right].
\]
Here, the right eigenvector---if the sum of priorities equals 100---is:
\[
\mathbf{w}^R \left( \mathbf{A} \right) = \left[
\begin{array}{K{2.5em} K{2.5em} K{2.5em} K{2.5em}}
    18.44 & 15.19 & 43.64 & 22.73 \\
\end{array}
\right].
\]
Consequently, the fourth alternative is ranked above the first.

The left eigenvector is
\[
\mathbf{w}^L \left( \mathbf{A} \right) = \left[
\begin{array}{K{2.5em} K{2.5em} K{2.5em} K{2.5em}}  
    24.82 & 38.78 & 10.49 & 25.91 \\
\end{array}
\right],
\]
hence, the elementwise reciprocal left eigenvector is
\[
\mathbf{w}^{-L} \left( \mathbf{A} \right) = \left[
\begin{array}{K{2.5em} K{2.5em} K{2.5em} K{2.5em}}
    20.14 & 12.89 & 47.67 & 19.29 \\
\end{array}
\right],
\]
implying that the first alternative is preferred to the fourth.

The consistency ratio of matrix $\mathbf{A}$ is $\mathit{CR}(\mathbf{A}) \approx 0.331$, which cannot be accepted according to Saaty's 10\% rule.

\citet{JohnsonBeineWang1979} have also generated 364 random matrices of order six with $\mathit{CI} < 1$ and found 164 ranking interchanges between the right and the componentwise reciprocal left eigenvectors. In addition, the authors note that the disagreement can occur for arbitrarily small positive values of $\mathit{CR}$ due to continuity.

For $n \geq 4$, the reciprocal property between the left and right eigenvector components holds not only if the matrix is consistent \citep{DeTurck1987}. In the case of the pairwise comparison matrix
\[
\mathbf{B} = \left[
\begin{array}{K{2em} K{2em} K{2em} K{2em}}
    1     &   8/5 &  1/4  & 4     \\
     5/8  & 1     &  5/8  & 10    \\
    4     &   8/5 & 1     & 4     \\
     1/4  & 1/10  &  1/4  & 1     \\
\end{array}
\right],
\]
the right eigenvector is 
\[
\mathbf{w}^R \left( \mathbf{B} \right) = \left[
\begin{array}{K{2em} K{2em} K{2em} K{2em}}
      2/9  &   5/18 &   4/9  &   1/18 \\
\end{array}
\right],
\]
and the left eigenvector is
\[
\mathbf{w}^L \left( \mathbf{B} \right) = \left[
\begin{array}{K{2em} K{2em} K{2em} K{2em}}
      1/4  &   1/5 &   1/8  &   1 \\
\end{array}
\right],
\]
thus, $\mathbf{w}^{-L} \left( \mathbf{A} \right) = \mathbf{w}^R \left( \mathbf{A} \right)$. However, $\mathit{CR}(\mathbf{B})$ is higher than $0.1$.

An acceptable inconsistency ratio does not guarantee that the right and inverse left eigenvectors imply the same ranking of the alternatives \citep{DoddDoneganMcMaster1995}.
Let
\[
\mathbf{C} = \left[
\begin{array}{K{2em} K{2em} K{2em} K{2em} K{2em}}
    1     & 1     & 3     & 9     & 9     \\
    1     & 1     & 5     & 8     & 5     \\
     1/3  &  1/5  & 1     & 9     & 5     \\
     1/9  &  1/8  &  1/9  & 1     & 1     \\
     1/9  &  1/5  &  1/5  & 1     & 1     \\
\end{array}
\right],
\]
where $\mathit{CR} \left( \mathbf{C} \right) \approx 0.082$.
The priorities from the right eigenvector are
\[
\mathbf{w}^R \left( \mathbf{C} \right) = \left[
\begin{array}{K{3em} K{3em} K{3em} K{3em} K{3em}}
    36.5652 & 38.9564 & 16.7155 & 3.4693 & 4.2936 \\
\end{array}
\right],
\]
and the priorities from the componentwise reciprocal of the left eigenvector are
\[
\mathbf{w}^{-L} \left( \mathbf{C} \right) = \left[
\begin{array}{K{3em} K{3em} K{3em} K{3em} K{3em}}
    40.6431 & 36.4208 & 15.0669 & 3.4391 & 4.4302 \\
\end{array}
\right].
\]
Therefore, the first alternative is the best by the left eigenvector but the second should be chosen by the right eigenvector.

\citet{IshizakaLusti2006} compare all the three weighting techniques---logarithmic least squares method, right, and inverse left eigenvectors---defined in Section~\ref{Sec2} for five intervals of the inconsistency ratio (0--0.02, 0.02--0.04, 0.04--0.06, 0.06--0.08, 0.08--0.1) and three to seven alternatives. The underlying matrices are generated randomly. Ranking contradictions are found to increase linearly with the level of inconsistency and the dimension of the matrix. However, the number of matrices in each case (100) is rather low, the matrix entries are integers, and the authors focus only on the ranking of the alternatives.

\citet[Section~6]{BozokiRapcsak2008} examine 100 million randomly generated pairwise comparison matrices of order five to estimate the frequency of rank reversal between the weights computed from the right and inverse left eigenvectors. This measure is found to increase along with the inconsistency ratio $\mathit{CR}$. However, the entries of the pairwise comparison matrix are chosen independently of each other, thus, the number of matrices for small values of $\mathit{CR}$ is relatively small. Furthermore, it remains to be seen whether the same result holds if the number of alternatives varies. 

\citet{Tomashevskii2015} argues that the derived weights cannot be used to rank the alternatives without taking their errors into account, and the rank reversal between the right and the inverse left eigenvectors occurs only because the errors are high. Consequently, the means do not contain any information on the ranking of the corresponding alternatives.

\section{The design of numerical experiments} \label{Sec5}

Our matrix generation algorithm is based on constructing a consistent pairwise comparison matrix and perturbing its elements. The technique of \citet{SzadoczkiBozokiJuhaszKadenkoTsyganok2022} is followed as it improves the method of \citet{SZadoczkiBozokiTekile2022}. This consists of the following steps:
\begin{enumerate}
\item
Choosing $n$ uniformly distributed random number $w_i$ from the interval $\left[ 1;9 \right]$. Computing a consistent pairwise comparison matrix $\mathbf{A} = \left[ a_{ij} \right]$, where $a_{ij} = w_i / w_j$ for all $i,j$.
\item
For all $i \neq j$, either $a_{ij}$ or $a_{ji}$ is perturbed depending on which entry is higher. \\
If $a_{ij} \geq 1$, then the perturbed entry $\hat{a}_{ij}$ is
\begin{equation} \label{Eq_perturbation}
\hat{a}_{ij} =
\begin{cases}
    a_{ij} + \varepsilon_{ij} & \text{if } a_{ij} + \varepsilon_{ij} \geq 1 \\
    1 / \left[ 1 - \varepsilon_{ij} - \left( a_{ij} - 1 \right) \right] & \text{otherwise},
\end{cases}
\end{equation}
where $\varepsilon_{ij}$ is a uniformly distributed random number from the interval $\left[ -\Delta; \Delta \right]$. \\
If $a_{ij} < 1$, then $a_{ji} > 1$, and $\hat{a}_{ji}$ is computed analogously to \eqref{Eq_perturbation}.
\item
The reciprocity of the matrix is kept by adjusting the pair of the perturbed entry.
\end{enumerate}
This process ensures that the perturbed elements are uniformly distributed around the original $a_{ij}$ on the scale where the distance between $1/b$ and $1/c$ equals the distance between $b$ and $c$, see \citet[Figure~1]{SzadoczkiBozokiJuhaszKadenkoTsyganok2022}.

We consider three different values of parameter $\Delta$ (1, 2, 3) and six dimensions ($4 \leq n \leq 9$) for which 1 million matrices are generated, respectively (altogether 18 million). For any matrix, the inconsistency ratio and the weight vectors according to the three techniques presented in Section~\ref{Sec2} are calculated.

\input{Figure2_simulated_matrix_CR_histogram}

The distribution of the pairwise comparison matrices with respect to the level of inconsistency is shown in Figure~\ref{Fig2}. If $\Delta = 1$, almost all matrices have an inconsistency ratio below $0.1$, except for $n = 4$. As the dimension of the matrix increases, the curves for a given $\Delta$ are more peaked and less asymmetric, meaning that inconsistency is more strongly determined by the maximal perturbation. This crucial observation is not reported in \citet{SzadoczkiBozokiJuhaszKadenkoTsyganok2022} since the authors provide only the average value of $\mathit{CR}$.
The reader is directed to \citet{SzadoczkiBozokiJuhaszKadenkoTsyganok2022} for other characteristics of the simulated matrices.

In order to compare two priority vectors $u$ and $v$ normalised by $\sum_{i=1}^n u_i = 1$ and $\sum_{i=1}^n v_i = 1$, respectively, four metrics are considered:
\begin{itemize}
\item
Euclidean distance:
\begin{equation} \label{Eq_M1}
d_{euc}(u,v) = \sqrt{\sum_{i=1}^n \left( u_i - v_i \right)^2};
\end{equation}
\item
Chebyshev distance:
\begin{equation} \label{Eq_M2}
d_{cheb}(u,v) = \max \left\{ \left| u_i - v_i \right|: 1 \leq i \leq n \right\};
\end{equation}
\item
Maximal ratio:
\begin{equation} \label{Eq_M3}
\omega(u,v) = \max \left\{ \max \left\{ \frac{u_i}{v_i}; \frac{v_i}{u_i} \right\}: 1 \leq i \leq n \right\};
\end{equation}
\item
Kendall tau:
\begin{equation} \label{Eq_M4}
\tau(u,v) = \frac{\# c(u,v) - \# d(u,v)}{n(n-1)/2};
\end{equation}
\end{itemize}
where $\# c$ ($\# d$) denotes the number of (dis)concordant pairs $1 \leq i,j \leq n$ when $u_i > u_j$ and $v_i > v_j$ ($v_i < v_j$).

The Euclidean distance \eqref{Eq_M1} is the length of a line segment between the two vectors.
The Chebyshev distance \eqref{Eq_M2} depends only on the greatest difference along any coordinate.

The maximal ratio \eqref{Eq_M3} is inspired by the Chebyshev distance but focuses on ratios instead of differences since the Chebyshev distance does not reflect high relative deviations in the weights of lower-ranked alternatives. For instance, let $u = \left[ 0.5;\, 0.4;\, 0.1 \right]$, $v = \left[ 0.5;\, 0.3;\, 0.2 \right]$, and $w = \left[ 0.6;\, 0.3;\, 0.1 \right]$. Then $d_{cheb}(u,v) = d_{cheb}(u,w) = 0.1$ but $\omega(u,v) = 2 > 4/3 = \omega(u,w)$. Clearly, the Chebyshev distance does not take into account whether the switching in priorities happens between the top two ($u$ and $w$) or the bottom two ($u$ and $v$) alternatives. However, the latter is more serious in terms of relative changes, which might be important if the priorities are used, for example, to allocate some resources proportionally.

Finally, Kendall tau \eqref{Eq_M4} is a standard rank correlation coefficient \citep{Kendall1938}. Its range is $\left[ -1;\, 1 \right]$ with $-1$ showing two opposite rankings and $+1$ indicating two identical rankings. In contrast to the other three measures, here a higher value is associated with a stronger agreement between the two weight vectors.

\section{Results} \label{Sec6}

For each of the four metrics \eqref{Eq_M1}--\eqref{Eq_M4} and each interval of the inconsistency ratio $\mathit{CR}$ (resolution: 0.005), we compute
\begin{itemize}
\item
the average value of the metric between the right eigenvector $\mathbf{w}^R$ and the inverse left eigenvector $\mathbf{w}^{-L}$ for all pairwise comparison matrices at the given level of inconsistency;
\item
the average value of the metric between the right eigenvector $\mathbf{w}^R$ and the (geometric) mean of the right and inverse left eigenvectors $\mathbf{w}^{\mathit{RL}}$ for all pairwise comparison matrices at the given level of inconsistency;
\item
the average value of the metric between the right eigenvector $\mathbf{w}^R$ and the row geometric mean weight vector $\mathbf{w}^{\mathit{RGM}}$ for all pairwise comparison matrices at the given level of inconsistency;
\item
the probability that the row geometric mean weight vector $\mathbf{w}^{\mathit{RGM}}$ is not farther from the right eigenvector $\mathbf{w}^R$ than the inverse left eigenvector $\mathbf{w}^{-L}$, based on all pairwise comparison matrices at the given level of inconsistency.
\end{itemize}
In order to ensure the robustness of the results, only the inconsistency intervals with at least one thousand matrices are shown; this accounts for the smaller range in the case of higher dimensions (cf.\ Figure~\ref{Fig2}).


\input{Figure3_average_Euclidean_distance}

\input{Figure4_average_Chebyshev_distance}

Figure~\ref{Fig3} uses the Euclidean, while Figure~\ref{Fig4} follows the Chebyshev approach to quantify the distances of the weights. In the case of average distances, the shapes of the lines are almost indistinguishable for a given $n$. For small levels of inconsistency, the distances between the right and inverse left eigenvectors increase almost linearly as a function of inconsistency, which coincides with the finding of \citet{IshizakaLusti2006}. However, the growth slackens and the maximum distance is reached for a given value of $\mathit{CR}$ somewhere between $0.22$ (if $n=9$) and $0.4$ (if $n=5$). The only exception is $n=4$, when a higher inconsistency is associated with a higher distance between the two eigenvectors.

On the other hand, the distance between the (right) eigenvector and the row geometric mean increases monotonically, albeit the curve is almost flat at higher levels of inconsistency. If $\mathit{CR}$ does not exceed 10\%, that is, inconsistency can be accepted according to the criterion of Saaty, then the logarithmic least squares method is essentially at the midpoint between the right and inverse left eigenvectors, and is close to the geometric mean of the two eigenvectors.
This seems to be a quite powerful argument for using the row geometric mean method to derive the weights.

Since the higher average distance between the two eigenvectors $\mathbf{w}^R$ and $\mathbf{w}^{-L}$ compared to the average distance between $\mathbf{w}^R$ and the row geometric mean $\mathbf{w}^{\mathit{RGM}}$ may be misleading, the probability that $\mathbf{w}^R$ is closer to $\mathbf{w}^{\mathit{RGM}}$ than to $\mathbf{w}^{-L}$ has also been calculated. If $\mathit{CR} < 0.1$, this is almost guaranteed for the Chebyshev distance, while the likelihood is still above 98\% for the Euclidean distance. The probability exceeds 80\% even at higher levels of inconsistency, and it is increasing with the number of alternatives at a given value of $\mathit{CR}$, except for $4 \leq n \leq 5$ and the Chebyshev distance.

\input{Figure5_average_maximal_ratio}

Figure~\ref{Fig5}, which focuses on the mean maximal ratios of the priorities, reinforces the findings from Figures~\ref{Fig3} and \ref{Fig4}. Consequently, the cardinal difference between the right and inverse left eigenvectors, as well as between the (right) eigenvector and the row geometric mean is robust and almost independent of the measure chosen.

\input{Figure6_average_Kendall_tau}

In contrast to the previous three metrics, Kendall tau considers only the rankings implied by the weights. Since a higher value shows a stronger similarity, Figure~\ref{Fig6} plots the difference of the average Kendall tau from its theoretical maximum of one. The number of ranking contradictions generally increases along with the level of inconsistency, but two breaks appear in the lines, especially for higher $n$, which are probably caused by the three different values for parameter $\Delta$. \citet[Fig.~3]{IshizakaLusti2006} present a similar result based on a much smaller number of randomly generated matrices.

However, according to \citet[Fig.~4]{IshizakaLusti2006}, the ranking contradiction phenomenon increases linearly with the dimension of the matrix because the possibility of a reversal rises if there are more alternatives. Figure~\ref{Fig6} does not fully support this conclusion: the value of mean Kendall tau is about the same for a given inconsistency interval if the number of alternatives is at least seven (e.g.\ about $0.9$ if $\mathit{CR}$ is approximately $0.1$). Again, the difference between the (right) eigenvector and the row geometric mean is only slightly larger than the difference between the right eigenvector and the (geometric mean) of the two eigenvectors. In addition, the ranking implied by the row geometric mean is not farther from the ranking implied by the right eigenvector than the ranking implied by the left eigenvector with a probability of at least 95\% (except for $n=4$, where this remains true only if $\mathit{CR} < 0.2$). Hence, the row geometric mean seems to be a reasonable compromise between the two eigenvectors even from an ordinal point of view.

Besides the results based on a high number of random pairwise comparison matrices, our simulations have provided some interesting examples that are worth further consideration.
First, according to \citet[p.~63]{JohnsonBeineWang1979}, the disagreement between the right and inverse left eigenvectors can occur at arbitrarily small positive values of inconsistency. For the pairwise comparison matrix
\[
\mathbf{M}^{(1)} = \left[
\begin{array}{K{3em} K{3em} K{3em} K{3em}}
    1     & 0.4759 & 0.9832 & 0.4025 \\
    2.1011 & 1     & 1.9975 & 0.7374 \\
    1.0171 & 0.5006 & 1     & 0.3704 \\
    2.4842 & 1.3560 & 2.6998 & 1 \\
\end{array}
\right],
\]
the weights from the right eigenvector are
\[
\mathbf{w}^R \left( \mathbf{M}^{(1)} \right) = \left[
\begin{array}{K{3em} K{3em} K{3em} K{3em} K{3em}}
    15.042 & 30.274 & 15.037 & 39.647 \\
\end{array}
\right],
\]
and the weights from the componentwise reciprocal of the left eigenvector are
\[
\mathbf{w}^{-L} \left( \mathbf{M}^{(1)} \right) = \left[
\begin{array}{K{3em} K{3em} K{3em} K{3em} K{3em}}
    15.036 & 30.281 & 15.049 & 39.635 \\
\end{array}
\right].
\]
Rank reversal arises between the first and the third alternatives in the last two positions, although $\mathit{CR} \left( \mathbf{M}^{(1)} \right) \approx 0.0007$.
Therefore, a class of pairwise comparison matrices with minimal inconsistency can be sought such that the right and inverse left eigenvectors lead to a different ranking of the alternatives, similar to the issue of Pareto inefficiency \citep{Bozoki2014}.

Second, the left and right eigenvector components might imply an opposite order of the alternatives. In the case of the pairwise comparison matrix
\[
\mathbf{M}^{(2)} = \left[
\begin{array}{K{3em} K{3em} K{3em} K{3em} K{3em}}
    1     & 1.624 & 0.574 & 1.072 & 1.054 \\
    0.616 & 1     & 1.132 & 1.089 & 1.269 \\
    1.743 & 0.884 & 1     & 1.515 & 0.467 \\
    0.933 & 0.919 & 0.660 & 1     & 1.694 \\
    0.949 & 0.788 & 2.140 & 0.590 & 1 \\
\end{array}
\right],
\]
the priorities from the right eigenvector are
\[
\mathbf{w}^R \left( \mathbf{M}^{(2)} \right) = \left[
\begin{array}{K{2.5em} K{2.5em} K{2.5em} K{2.5em} K{2.5em}}
    19.75 & 19.16 & 20.85 & 19.53 & 20.71 \\
\end{array}
\right],
\]
and the priorities from the componentwise reciprocal of the left eigenvector are
\[
\mathbf{w}^{-L} \left( \mathbf{M}^{(2)} \right) = \left[
\begin{array}{K{2.5em} K{2.5em} K{2.5em} K{2.5em} K{2.5em}}
    20.25 & 20.55 & 19.31 & 20.27 & 19.62 \\
\end{array}
\right].
\]
Hence, the ranking of the alternatives is $3 \succ 5 \succ 1 \succ 4 \succ 2$ in the former case, which is reversed to $2 \succ 4 \succ 1 \succ 5 \succ 3$ in the latter case.
The inconsistency ratio is $\mathit{CR} \left( \mathbf{M}^{(2)} \right) \approx 0.078$.

Third, it has been thought that ``\emph{Only very close priorities suffer from ranking contradictions}'' \citep[p.~398]{IshizakaLusti2006}. Even though the exact meaning of very close remains obscure, the following example probably disproves this statement:
\[
\mathbf{M}^{(3)} = \left[
\begin{array}{K{3em} K{3em} K{3em} K{3em} K{3em}}
    1     & 0.371 & 2.013 & 5.389 & 0.243 \\
    2.698 & 1     & 4.596 & 7.527 & 0.736 \\
    0.497 & 0.218 & 1     & 2.321 & 0.167 \\
    0.186 & 0.133 & 0.431 & 1     & 0.385 \\
    4.120 & 1.359 & 5.973 & 2.598 & 1 \\
\end{array}
\right],
\]
where the weights from the right eigenvector are
\[
\mathbf{w}^R \left( \mathbf{M}^{(2)} \right) = \left[
\begin{array}{K{2.5em} K{2.5em} K{2.5em} K{2.5em} K{2.5em}}
    15.26 & 33.23 & 7.74  & 5.68  & 38.08 \\
\end{array}
\right],
\]
and the weights from the componentwise reciprocal of the left eigenvector are
\[
\mathbf{w}^{-L} \left( \mathbf{M}^{(2)} \right) = \left[
\begin{array}{K{2.5em} K{2.5em} K{2.5em} K{2.5em} K{2.5em}}
    15.29 & 37.84 & 8.55  & 4.93  & 33.39 \\
\end{array}
\right].
\]
Consequently, the best alternative is the fifth according to the right and the second according to the inverse left eigenvector. The absolute value of the difference between the weights of the second and the fifth alternatives is $4.85$ and $4.44$, respectively, if the sum of weights is normalised to 100.
The inconsistency ratio is $\mathit{CR} \left( \mathbf{M}^{(3)} \right) \approx 0.0993$.

\section{Conclusions} \label{Sec7}

The paper has addressed one of the most serious shortcomings of the eigenvector method, a widely used priority deriving technique for pairwise comparison matrices. In particular, we have compared the weights implied by the right and inverse left eigenvectors, as well as by the row geometric mean for a large set of randomly generated matrices. The two eigenvectors turned out to lead to different priorities and rankings relatively often, and their disagreement is not necessarily a monotonic function of the level of inconsistency and the number of alternatives. Although the problem is less threatening if inconsistency remains below the acceptable threshold, the left and right eigenvector components may imply an opposite priority order or a rank reversal between alternatives with distant weights even if the value of the inconsistency ratio does not exceed 10\%.

There are at least two important conclusions to be drawn from these results.
Our findings uncover that the row geometric mean is almost at the midpoint between the principal right and inverse left eigenvectors, hence, it is a reasonable compromise between them. Therefore, a novel argument is provided for following this method, which is also easy to calculate and satisfies many attractive theoretical properties.
Furthermore, we have reinforced an important message of \citet{JohnsonBeineWang1979}: assuming the existence of a single ranking or priority vector is usually too demanding. Rather, a range of possible orderings and weights can be allowed (instead of a single ranking and exact priorities) such that the range is wider for a higher level of inconsistency and uncertainty in the input data.

Finally, this research is far from finished and can be continued in several directions.
First, the disagreement between the two eigenvectors is worth analysing on particular classes of pairwise comparison matrices. For instance, following the work of \citet{DeTurck1987}, it would be interesting to characterise the set of matrices for which the reciprocal property between the left and right eigenvector components hold.
Second, a new weighting method can be introduced by aggregating the two eigenvectors appropriately.
Third, other deficiencies of the eigenvector method, such as Pareto inefficiency \citep{BlanqueroCarrizosaConde2006, BozokiFulop2018}, might be studied in a similar Monte Carlo experiment. 

\section*{Acknowledgements}
\addcontentsline{toc}{section}{Acknowledgements}
\noindent
We are grateful to \emph{S\'andor Boz\'oki} and \emph{Zsombor Sz\'adoczki} for useful advice. \\
Three anonymous reviewers provided valuable comments and suggestions on earlier drafts. 

\bibliographystyle{apalike}
\bibliography{All_references}

\end{document}